\documentclass{jnmp01b}

\usepackage{amsmath}

\numberwithin{equation}{section}

\theoremstyle{plain}
  \newtheorem*{conditionA}{Condition A}
  \newtheorem*{conditionB}{Condition B}
  \newtheorem*{conditionC}{Condition C}
  \newtheorem*{definition}{Definition}
  \newtheorem{lemma}{Lemma}
  \newtheorem{theorem}{Theorem}
\theoremstyle{definition}
  \newtheorem*{remark}{Remark}
  \newtheorem{example}{Example}

\newcommand{\D}{\mathbb}
\newcommand{\DS}{\displaystyle}
\def \c{\cite}
\def \r{\ref}
\def \l{\label}

\makeatletter
\DeclareRobustCommand{\primfrac}[1]{%
   \PackageWarning{amsmath}{%
Foreign command \@backslashchar#1; %
\protect\frac\space or \protect\genfrac\space should be used instead%
   }
   \global\@xp\let\csname#1\@xp\endcsname\csname @@#1\endcsname
   \csname#1\endcsname
}
\makeatother

\begin{document}

\renewcommand{\evenhead}{I.\ Anders and A.\ Boutet de Monvel}
\renewcommand{\oddhead}{Asymptotic Solitons of the  Johnson Equation}


\thispagestyle{empty}

\begin{flushleft}
\footnotesize \sf
Journal of Nonlinear Mathematical Physics \qquad 2000, V.7, N~3,
\pageref{firstpage}--\pageref{lastpage}.
\hfill {\sc Article}
\end{flushleft}

\vspace{-5mm}

\copyrightnote{2000}{I.\ Anders and A.\ Boutet de Monvel}

\Name{Asymptotic Solitons of the  Johnson Equation}

\label{firstpage}

\Author{Igor ANDERS~$^\dagger$
               and Anne BOUTET de MONVEL~$^\ddagger$}

\Adress{$^\dag$ Mathematical Division,
         Institute for Low Temperature Physics,
         47 Lenin Avenue, \\
         ~~310164 Kharkov, Ukraine\\[2mm]
         $^\ddag$ Universit\'e Paris-7,
         Physique math\'ematique et G\'eom\'etrie,
         Institut de Math\'ematiques, \\
         ~~case 7012,
         2 place Jussieu, 75251 Paris Cedex 05, France}

\Date{Received November 25, 1999; Accepted February 29, 2000}

\begin{abstract}
\noindent
We prove the existence  of   non-decaying real solutions of
the Johnson
equation, vanishing as $x\to+\infty$.
We obtain  asymptotic formulas as $t\to\infty$  for the solutions in the
form of an infinite series of
asymptotic solitons with curved lines of constant phase and varying
amplitude and width.
\end{abstract}




\section{Introduction}

The Johnson equation (JE)
\begin{equation}
\left (v_t+{1\over 4}v_{xxx}+{3\over 2}vv_x +{v\over 2t}\right
)_x=-{12\alpha^2\over
   t^2}u_{yy}
\l{ej}
\end{equation}
$(\alpha^2=\pm 1)$ or the cylindrical Kadomtsev-Petviashvili equation, is the
  analogue of the well-known cylindrical Korteweg-de Vries
equation $(\alpha =0)$ in two spatial dimensions (2D). The JE was
obtained firstly in \c{J} under the description of the surface waves on
a shallow incompressible liquid. Later it was shown that it describes the
propagation of
waves in the stratified media \c{L}. It follows from the derivation of
(\r{ej}) that the correct statement of the Cauchy problem is possible
only as $t=t_0>0$.

In \c{LMS} and \c{MS} the equivalence of the Kadomtsev-Petviashvili
equation (KP) and the JE was established.
Let $u(\xi,\eta,\tau)$ be an arbitrary  solution of the KP
\begin{equation}
\left (u_\tau +{1\over 4}u_{\xi\xi\xi}+{3\over 2}uu_\xi\right )_\xi
=-{3\alpha^2\over 4}u_{\eta\eta}.
\l{kp}
\end{equation}
Then the function
\begin{equation}
v(x,y,t)=u\left (x -{y^2t\over 48\alpha^2},{yt\over 4},t \right )
\l{v}
\end{equation}
satisfies the JE. This mapping $u(\xi,\eta,\tau)\to v(x,y,t)$ is
invertible. Each solution $v(x,y,t)$ of the JE generates a solution
of the KP by the formula
\begin{equation}{}
u(\xi,\eta,\tau)=v\left(\xi+{\eta^2\over 3\alpha^2 t}, {4\eta\over
\tau},\tau\right ).
\l{u}
\end{equation}
It was shown in \c{LMS} that mappings (\r{v}) and (\r{u}) preserve the
class of functions rapidly decaying at infinity (as $(x^2+y^2)^{-1}\to
0$), and all the results obtained in the theory of the KP in the
corresponding class of solutions can be directly applied to solve the
Johnson equation. The situation is similar for the solutions
involving the Airy functions (see \c{SM} and \c{N}).
Obviously this is not the fact  for periodic
initial data, which are not
invariant with respect to this transformation, and investigation of
the JE is an independent interesting problem in this case
(for the KP see the corresponding theory, for example, in
\c{kri}-\c{m3}).

We are interested in the construction of
a class of JE-I ($\alpha =i$ in (\r{ej})) non-decaying solutions,
which are bounded for all
$(x,y,t)$ and vanish as $x\to +\infty$ for all fixed $y$ and $t$.
Such a kind of KP solutions was constructed and investigated
firstly for KP-II in \c{akk}, \c{A3}, and then for KP-I in
\c{A1}-\c{A2}. It turns out that all basic  stages of the construction of the
solutions of the KP and the JE, and the study of their asymptotic
behaviour admit mutual recounting using the described mapping.

We apply the change of variables
\[
\xi= x -{y^2t\over 48\alpha^2},\quad \eta ={yt\over 4},\quad \tau =t
\]
to the scheme of the
V.E.Zakharov and A.B.Shabat ``dressing method''  \c{ZS} of integration
of the KP, and obtain analogous formulas
for the JE. Using them we prove the existence of a class of JE-I
non-decaying
solutions with the prescribed properties. The simplest one is the
one-soliton solution
\[
v(x,y,t)={2q^2\over \cosh^2\left [q\left (x-\left
(q^2-3p^2-{y^2\over 48}-{py\over 2}\right )t-{1\over 2q}
\ln{c\over 2q} \right )\right ]}
\]
$(p\in\D{R}$, $q\in\D{R}^+)$ which corresponds to the KP
plane-soliton by virtue of (\r{v}).

Then we study the asymptotic behaviour of the constructed solution as
$t\to\infty$.
The investigations of long-time asymptotic behaviour of non-decaying
solutions of  2D non-linear
evolution equations
is closely connected with the same investigations in one spatial dimension.
A.V.Gurevich and
L.P.Pitaevsky studied  in 1973 a non-decaying  solution of the
Korteweg--de Vries (KdV)
equation, which describes the evolution of an initial step--function
(\c{GP1},\c{GP2}).
They applied Whitham method  to construct an approximation
of this solution by a  knoidal wave with
slowly varying parameters and detected   the appearance of
many strong oscillations like solitons on the front of
the solution for a large time.
This approximate solution satisfies the KdV--equation with error
vanishing as $t\to +\infty$.
The mathematical ground of this phenomenon was done in 1975 by
E.Ya.Khruslov in \c{Kh1} and \c{Kh2},
where the nature of these solitons was explained.  Subsequently these
solitons were called asymptotic solitons.
An analogous phenomenon of splitting of non-decaying initial data into
infinite series of
solitons was proved later for other KdV--like equations  (nonlinear
Schr{\"o}dinger equation,
sine--Gordon equation, modified KdV and the Toda lattice as a discrete
analogue of the KdV) (\c{KK}--\c{ABE}).

In \c{akk}--\c{A2}  the method proposed by E.Ya.Khruslov  (\c{Kh1},\c{Kh2})
was extended to the investigation
of the asymptotic behaviour of non-decaying solutions of KP-type
equations (KP, modified KP-I and 2D-Gardner equation) as
$t\to\infty$.  It was proved that they are represented as infinite series
of solitons with
curved lines of constant phase in the neighbourhood of the front as
$t\to\infty$. These
asymptotic solitons were called curved asymptotic solitons. Note
that
recently (\c{z}) V.E.Zakharov
also considered a curved soliton of the KP-II equation, but in another
space--time
domain.

Our principal goal is to prove the phenomenon of splitting of non-decaying
solutions of the JE-I into infinite series of solitons as $t\to\infty$.

We prove that there exist non-decaying real solutions of
the JE-I, which split in the neighbourhood of the front into a series of
solitons of the form
\begin{equation}
v_n(x,y,t)=
\displaystyle{{2q_0(y)^2\over
     \cosh^2 \left [q_0(y)\left ( x-C(y)t+{1\over 2q_0(y)}\left (\ln
t^{n+1/2}-\ln g(y)-
   \phi_n(y)\right )          \right ) \right]}},
\l{vn}
\end{equation}
which depend on two parameters $C(y)$ and $g(y)$. The functions
$p_0(y)$, $q_0(y)$ and $\phi_n(y)$  are completely determined by
them. These solitons are diverged with the velocity
$\ln t^{{n+1/2\over 2q_0(y)}}$. They have  varying amplitude and width
in the general case, but we present also  examples where amplitude and
width  are
constant. In these cases curved and weakly curved
asymptotic solitons are both constructed. The lines of constant phase of
the weakly curved solitons are deviated from the straight line just on a
value proportional to $\ln y^2$.

Asymptotic solitons (\r{vn}) of the JE-I and the KP-I (\c{A1}, \c{A12})
coincide taking into account transformations (\r{v}) and (\r{u}).

\section{Existence of a Johnson equation solution.}

After application of  the change of variables
\[
\xi= x -{y^2t\over 48\alpha^2},\quad \eta ={yt\over 4},\quad \tau =t
\]
to the scheme of the V.E.Zakharov and A.B.Shabat ``dressing method''  \c{ZS}
of integration of the KP we  obtain the following formulas
for the JE.
A JE solution has the form
\begin{equation}
v(x,y,t)=2{\partial\over \partial x} K(x,x,y,t),
\l {1.1}
\end{equation}
where the function $K(x,s,y,t)$ is a solution of the Marchenko integral
equation
\begin{equation}
K(x, z, y, t)+F(x, z, y, t)+\int_x^\infty K(x, \xi , y, t)F(\xi, z, y, t)d\xi
=0.
\label{1.2}
\end{equation}
This equation is an equation with respect to $z$ and $x,y,t$ are parameters.
The kernel  $F(x,z,y,t)$ of (\r{1.2}) satisfies the system of linear
differential equations
\begin{equation}
\begin{cases}
\DS
F_{t} +{y^2\over 48\alpha^2}(F_x+F_z)+{y\over
     4\alpha}(F_{xx}-F_{zz})+  F_{xxx} + F_{zzz}=0\\[1em]
\DS
\alpha F_y+{yt\over
     24\alpha}(F_x+F_z)+ {t\over 4}\left (F_{xx}-F_{zz}\right )=0
\end{cases}
\l {1.3}
\end{equation}
($\alpha =i$ for the JE-I, and $\alpha =1$ for the JE-II).
The correspondence of (\r{1.1})-(\r{1.3}) to the JE can be verified
by a direct method, described in \c{AS} for the KP.
The scheme (\r{1.1})-(\r{1.3})  don't allow us to solve the Cauchy
  problem for JE, but
it is rather convenient for the construction of classes of solutions
with various properties, in particular, of rapidly decaying
rational solutions described in \c{SM} and of non-decaying
solutions.

A wide class of solutions of (\r{1.3}) as $\alpha =i$ can be found
as follows
\begin{equation}
F(x,z,y,t)
=\iint_{\Omega} \exp[ip(x -z )-q(x +z )+2qf(p,q,y)t]d\mu(p,q),
\l {1.4}
\end{equation}
where $\Omega\subset\D{C}^+$ ($\D{C}^+=\{\lambda\mid\lambda =
  p+iq, q>0\}$ is the upper half-plane of the complex plane),
$f(p,q,y)=q^2-3p^2+py/2-y^2/48$,
and  $d\mu(p,q)$ is some measure on $\Omega$.

To construct a JE-I solution by the scheme (\r{1.1})-(\r{1.4}) we must
define
the set $\Omega$ in (\r{1.4}) and the measure $d\mu(p,q)$ on this
set. For this goal  we  introduce
two  positive functions $C(s)$ and $g(s)$ which play an important role
in the construction of the solution and investigation of its
asymptotic behaviour, and we formulate the following conditions.

\begin{conditionA}
The function $C(s):\D{R}\to\D{R}^+$ is of class $C^2$
and such that
\begin{equation}
C(s)\geq\delta >\varepsilon^2\quad(\delta,\varepsilon =\mbox{\rm
const}>0),\quad C^{\prime\prime}(s)>-1/24.
\l{1.5}
\end{equation}
\end{conditionA}

\begin{conditionB}
The set $\Omega$ has the form
\begin{equation}
\Omega =\left\{(p,q)\in\D{R}^2\mid -\infty <p<\infty,\; 0<\varepsilon \leq
                q\leq h(p)\right \},
\l {1.7}
\end{equation}
where $q=h(p)$ is the envelope of the family of hyperbolas
\begin{equation}
f(p,q,s)=C(s),
\l {1.6}
\end{equation}
which touch it at the point
\[
\left (p_0(y),q_0(y)\right )=\left (C^{\prime}
(y)+y/12,\,\,
\sqrt{C(y)+3\left (C^\prime (y)\right )^2}\right ).
\]
\end{conditionB}

\begin{remark}
The special structure of $\Omega$ (\r{1.7})
implies that
\begin{equation}
C(s)=\max_{(p,q)\in\Omega}f(p,q,s).
\l{max}
\end{equation}
\end{remark}

\begin{conditionC}
The function $g(s):\D{R}\to\D{R}^+$, $g(s)<A=\mbox{\rm const}$
is of class $C^\infty$ and
such that the measure $d\mu$ of the form $d\mu(p,q)=\tilde g(p,q)dpdq$ with
real positive $\tilde g\in C^\infty$,  $\tilde g(p_0(Y),q_0(Y))=g(Y)$,
satisfies the inequality
\begin{equation}
\forall a=\mbox{\rm const}>0\,:\quad
\iint_{\Omega}e^{a(q+|p|)}d\mu(p,q) <\infty.
\l{1.8}
\end{equation}
\end{conditionC}

Let us show that under Conditions A-C the scheme (\r{1.1})-(\r{1.3})
determines a smooth real
solution  of the JE-I  vanishing as $x\to+\infty$. For a function
$h(y)\in L^2[x,\infty)$ define the operator $\hat F$ by
\begin{equation}
[\hat Fh](z)=\int_x^{\infty} F(s,z,y,t)h(s)ds
\l{fh}
\end{equation}
with the kernel $F$ given by (\r{1.4}), where  $h(s)$ also depends on
the parameters $y,t$.

\begin{lemma}
Assume that Conditions A-C are fulfilled.\\
Then  $\hat F$ is a self-adjoint, compact and positive operator in
$L^2[x,\infty)$.
\end{lemma}

\begin{proof}
Its self-adjointness follows from
the form (\r{1.4}) of $F(x,z,y,t)$.
Let us show that this
operator is compact. We estimate the Hilbert-Schmidt norm of  $\hat F$:
\begin{align*}
||\hat F||^2_{L^2[x,\infty)}
&=\int_x^{\infty}\int_x^{\infty}|F(x,z,y,t)|^2dsdz\\
&\leq
{1\over 4\varepsilon^2}\iint_{\Omega}d\mu (p,q)
\iint_{\Omega}\exp
[4q(|x|+f(p,q,y))t]d\mu (p,q)\\
&\leq {1\over 4\varepsilon^2}\iint_{\Omega}d\mu (p,q)
\iint_{\Omega}e^{aq}d\mu (p,q)< \infty.
\end{align*}
Thus $\hat F$ is a
Hilbert-Schmidt operator. Hence, it is a compact operator (\c{KF}).
To prove the positivity of $\hat F$, consider the scalar product
\begin{align*}
(\hat Fh,h)
&=\int_x^\infty\int_x^\infty
F(s,z,y,t)h(s)ds\overline h(z)dz\\
&
=
\iint_{\Omega}e^{2qf(p,q,y)t}\biggl| \int_x^\infty
e^{(ip-q)s}h(s)ds\biggr|^2d\mu(p,q) >0
\end{align*}
as $h(s)\ne 0$.
\end{proof}

Under the conditions of Lemma 1 the following statement holds.

\begin{lemma}
The scheme $(\r{1.1})$-$(\r{1.4})$ determines a smooth real solution of
the JE-I vanishing as $x\to\infty$ and bounded for all
fixed $x,y,t$ $(t>0)$.
\end{lemma}

\begin{proof}
Let us represent (\r{1.3}) in the operator form in $L^2[x,\infty)$
\begin{equation}
\varphi+\hat F\varphi=f,
\l{ie}
\end{equation}
where $\hat F$ has the form (\r{fh}), and $\varphi =K(x,z,y,t)$,
$f=-F(x,z,y,t)$.
Due to the positivity of $\hat F$ the homogeneous equation $\varphi
+\hat F\varphi =0$ has only the trivial solution. Since $\hat F$ is a
compact operator, then by the Fredholm theorem (\c{KF}) inhomogeneous
equation (\r{ie}) has a unique solution given by
\[
K(x,z,y,t)=-(I+\hat F)^{-1}F(x,z,y,t)
\]
with $\vert |(I+\hat F)^{-1}\vert |\leq 1$.

Due to Condition C and the fact that $\Omega$ is inside
the upper half-plane at positive
distance from the $q$-axis (Condition B),
$F$ is an infinitely differentiable function
with respect to all variables. Moreover, $D_i^\alpha F\to 0$
$(D^\alpha ={\partial^\alpha \over \partial x_i^\alpha };\,\, \alpha
=0,1,\dots;\,\, i=1,\dots,4)$ as
$x+z\to\infty$, and  $D_i^\alpha F$ are bounded for all fixed  $x,z,y,t$
$(t>0)$. One can show (\c{VA})
that the function $K$ has the same properties.

Let us prove that $K(x,x,y,t)$ is a real function. After
multiplication of (\r{ie}) by $\overline \varphi$ and integration with
respect to $z$ from $x$ to $+\infty$ we obtain
\begin{equation}
\vert | \varphi\vert|^2+(\hat F\varphi,\varphi)=(f,\varphi).
\l{iih}
\end{equation}
The self-adjointness of $\hat F$ implies that the imaginary part of
the left-hand-side of (\r{iih}) is equal to zero:
\begin{equation}
{\rm Im}\int_x^\infty F(x,s)\overline{K(x,s)}ds =0.
\l{imi}
\end{equation}
Application of conjugation to (\r{1.3}) for $z=x$ gives
\[
\overline{K(x, x)}+F(x, x)+\int_x^\infty \overline{K(x,
\xi)}F(x,\xi)d\xi =0.
\]
It follows from (\r{imi}) and the reality of $F(x,x,y,t)$ that
$K(x,x,y,t)$ is real.
\end{proof}

\section{Theorem about long-time asymptotic behaviour
          of JE-I non-decaying solutions}

Our goal is to investigate the long-time asymptotic behaviour of the JE-I
solution defined in the previous section. To define a domain in which
we shall carry out the investigations,
we introduce the following definition.

\begin{definition}
Let $M>2$ be an arbitrary number. The domain
$G_M(t)\subset\D{R}^2$ given by
\[
G_M(t)=\left\{(x,y)\in \D{R}^2\,\Bigm|\,
|\ln g(y)|<\ln t, \,\,\,
x>C(y)t-{1\over 2q_0(y)}\ln{ t^{M+1}}\right\}
\]
is called the  neighbourhood of the solution front.
\end{definition}

The following theorem describes the asymptotic
behaviour of the JE-I solution defined by Lemma 2 for
large time.

\begin{theorem}
Assume that  Conditions A-C are fulfilled.

\noindent
Then the JE-I solution $v(x,y,t)$ constructed by the scheme
{\rm(\r{1.1})--(\r{1.4})}
is represented in the domain $G_M(t)$  as $t\to\infty$ in the
following way
\begin{gather}
v(x,y,t)=\sum_{n=1}^{[M-1]}v_n(x,y,t)+
O\left ({1\over t^{1/2-\varepsilon_1}}\right ),\qquad
(0<\varepsilon_1 <1/2)
  \l{3.3}
\\
v_n(x,y,t)=
\displaystyle{{2q_0(Y)^2\over \cosh^2\left [q_0(Y)
\left (x-C(Y)t+{1\over 2q_0(Y)}\left (\ln t^{n+1/2}-\ln g(Y)-
\ln \phi_n(y)
\right )\right)\right]}},\notag
\end{gather}
where $q_0(Y)=\sqrt{C(y)+3\left (C^\prime (y)\right)^2}$,
\[
\phi_n(y)={(C(y)+48(C^\prime(y))^2)^{n-1}(1+24C^{\prime\prime}(y))^{n-1/2}
Q^{(n)}\Gamma^{(n)}\over
2^{(2n+5)/2}((n-1)!)^2(C(y)
+12(C^\prime(y))^2)^{(10n-3)/4}Q^{(n-1)}\Gamma^{(n-1)}},
\]
and $\Gamma^{(n)},Q^{(n)}>0$ are the determinants of the $n$ by $n$ matrices
with entries
\begin{gather*}
\Gamma^{(n)}_{i+1,k+1}=\Gamma \left( {i+k+1\over 2}\right)(1+(-1)^{i+k}),
\quad
Q^{(n)}_{i+1,k+1}=\Gamma (i+k+1),\\
i,k=0,\dots,n-1.
\end{gather*}
Here the asymptotic representation $(\r{3.3})$ is uniform with respect to $x$
and $y$ in $G_M(t)$ for any fixed $M\geq 2$.
\end{theorem}

Let us mark the key points of the proof.
The proof consists in three steps.
On the first step we show that as $t\to\infty $
the kernel $F(x,z,y,t)$ of integral equation (\r{1.2}) is represented
as the sum of a degenerate kernel and a kernel with small operator norm
in the
space $L^2[x,\infty)$.
On the second step we prove that the degenerate kernel
gives the main contribution in the
asymptotic representation  of the  solution of equation (\r{1.2}).
The third step consists in the analysis of representation (\r{1.1})
for the solution $v(x,y,t)$ as $t\to\infty$, where the function
$K(x,z,y,t)$
is a solution of the Marchenko integral equation with
the degenerate kernel.

\bigskip

\noindent
\underline{First step.}
To investigate the Marchenko equation kernel $F(x,z,y,t)$ (\r{1.4}) as
$t \rightarrow \infty$, we set
$x=C(y)t+\xi, \quad z=C(y)t+\zeta $
and define $\tilde F(\xi ,\zeta ,y,t)=F(\xi +C(y)t ,\zeta +C(y)t,y,t)$.
Then the function
$\tilde F(\xi ,\zeta ,y,t) $ is written as follows:
\[
\tilde F(\xi ,\zeta ,y,t)=
\iint_{\Omega}\exp[ip(\xi -\zeta )-q(\xi +\zeta )-
2q(C(y)-f(p,q,y))t]d\mu(p,q).
\]
For sufficiently small $\varepsilon^\prime >0$ let us consider the curve
\begin{equation}
2q(f(p,q,y)-C(y))+\varepsilon^\prime =0.
\l{3.4}
\end{equation}
This curve separates the domain $\Omega$ into two subdomains
$O_{\varepsilon^\prime}$
and $\Omega_{\varepsilon^\prime}$ so that
$\Omega=\overline{O_{\varepsilon^\prime}}
\cup \Omega_{\varepsilon^\prime}$. Here  ${O_{\varepsilon^\prime}}$
lies between the curve $q=h(p)$ and curve (\r{3.4}),
moreover $(p_0(y),q_0(y)) \in O_{\varepsilon^\prime}$.
The set $\Omega_{\varepsilon^\prime}$ is the
complement of the set $\overline{O_{\varepsilon^\prime}}$
in the domain $\Omega$. According to this decomposition,
kernel (\r{1.4}) is the sum of two kernels
which we denote $F_1(x,z,y,t)$ and $F_2(x,z,y,t)$
respectively.

Let us make a change of variables, setting
\begin{equation}
r=2q(C(y)-f(p,q,y))
\l{3.5}
\end{equation}
in the kernel $F_1(x,z,y,t)$ which contains integration  over the set
$O_{\varepsilon^\prime}$. Let $u$ be the projection of a radius vector
directed
from the point
$(p_0(y),q_0(y))$ to the point $(p,q)\in O_{\varepsilon^\prime}$ on the
tangent
to the curve
$h(p,q)=0$ at the point $(p_0(y),q_0(y))$ or, that is the same,
on the tangent to the curve $f(p,q,y)=C(y)$ at the same point, i.e.:
\begin{equation}
u={12p_0(y)-y\over\sqrt{16q_0^2(y)+(12p_0(y)-y)^2}} (q-q_0(y)) +
{4q_0(y)\over\sqrt{16q_0^2(y)+(12p_0(y)-y)^2}} (p-p_0(y)).
\l{3.6}
\end{equation}

The system of equations (\r{3.5}), (\r{3.6})
has a unique solution with respect to $p$ and $q$ in $O_\varepsilon$
as $\varepsilon^\prime\leq {2\delta^3\over\sqrt 3}$.
Therefore in the neighbourhood of the point $(p_0,q_0)$
the variables $p$ and $q$ can be expressed via the variables $r$ and $u$:
\begin{align*}
p(r,u)&=p_0+k_1r+k_2u+k_3ur+k_4r^2+k_5u^2\dots,\\
q(r,u)&=q_0+\lambda_1r+\lambda_2u+\lambda_3ur+\lambda_4r^2+\lambda_5u^2\dots,
\end{align*}
where $k_n$, $\lambda_n$ are the coefficients of the corresponding Taylor
series. The first have the form
\begin{align*}
k_1(y)&={ 12p_0(y)-y\over q_0(y)(16q_0^2(y)+(12p_0(y)-y)^2)},\\
k_2(y)&={4q_0(y)\over\sqrt{16q_0^2(y)+(12p_0(y)-y)^2}},\\
\lambda_1(y)&=-{4\over 16q_0^2(y)+(12p_0(y)-y)^2},\\
\lambda_2(y)&={ 12p_0(y)-y\over\sqrt{16q_0^2(y)+(12p_0(y)-y)^2}}.
\end{align*}
One can obtain the expansion coefficients $k_n$ and $\lambda_n$ in an
explicit form after $n$--times differentiation of (\r{3.5}) and (\r{3.6})
with respect to $r$ and $u$.
In the neighbourhood of the point
$(p_0(y),q_0(y))$ the equation $q=h(p)$ can be written using variables
$r$ and $u$. It is easy to check that
${\partial^2\over\partial p^2}h(p)\vert_{p=p_0(y)}\ne 0$
since $C^{\prime\prime}> -1/24$ (Condition A). Therefore the curves
$f(p,q,y)=C(y)$ and $q=h(p)$ have a contact of
the first order, and $q=h(p)$ takes the form $u=u(r)$:
\[
u=\pm a(y)\sqrt r+b(y)r+\dots,
\]
where
\begin{gather*}
a(y)=\left [ {16q_0^2(y)+(12p_0(y)-y)^2\over
2q_0(y)\left (48q_0^2(y)-(12p_0(y)-y)^2-16h_{pp}(p_0(y))q_0^3(y)\right )}
\right ]^{1\over 2},\\
h_{pp}(p_0)={\partial^2 h(p)\over\partial p^2}\vert_{p=p_0}.
\end{gather*}
Using the new variables $u$ and $r$
and the notation
$E_0(\xi,\zeta,y)= e^{ip_0(y)(\xi-\zeta)-q_0(y)(\xi+\zeta)}$
we write the function
$\tilde F_1(\xi,\zeta,y,t)=F_1(\xi+C(y)t,\zeta+C(y)t,y,t)$
as follows:
\begin{gather}
\tilde F_1(\xi,\zeta,y,t)=E_0(\xi,\zeta,y) \notag\\
\qquad\qquad\quad\times\int^{\varepsilon^\prime}_{0}dr
\int^{a\sqrt r +\dots}_{-a\sqrt r+\dots}
\,du\, j(r,u,y)\tilde g(r,u,y)
e^{i(p-p_0)(\xi-\zeta)-(q-q_0)(\xi +\zeta) - rt}\l{3.7}
\end{gather}
where $j(r,u,y)=j(p(r,u,y),q(r,u,y))$ is the Jacobian corresponding to
the change of variables $(p,q)\to (r,u)$.
Let us expand integrand in (\r{3.7}) into a series with respect to the
powers of $r$ and
$u$ in the neighbourhood of the point
$(p_0(y),q_0(y))$ $(u=0,r=0)$:
\begin{gather*}
j(r,u,y)\tilde g(r,u,y)\exp[i(p-p_0(y))(\xi-\zeta)-(q-q_0(y))(\xi +\zeta )]
\\
\qquad=\sum_{n=0}^{\infty}
\sum_{j=0}^{n}\sum_{l=0}^{j}\sum_{m=0}^{n-j}\zeta^j\xi^{n-j} r^{l+m}
u^{n-l-m}\varphi_{n,j,l,m}(y) (1+\psi_n(r,u)),
\end{gather*}
where
\begin{gather*}
\varphi_{n,j,l,m}(y) = {(-1)^{n-m} \over l!m!(n-j-m)!(j-l)!}
(ik_1(y)+\lambda_1(y))^l(ik_2(y)+\lambda_2(y))^{j-l} \\
\qquad\qquad\qquad\times(ik_1(y)-\lambda_1(y))^m
(\lambda_2(y)-ik_2(y))^{n-j-m} g(y)j_0(y), \\
j_0(y)=j(0,0,y)={1\over q_0(y)\sqrt{16q_0^2(y)+(12p_0(y)-y)^2}},
\end{gather*}
$g(y)=\tilde g(0,0,y),$ and the functions
$\psi_n(r,u)$ satisfy $|\psi_n(r,u)|\leq An(r+|u|)$.

After integration with respect to $u$ and $r$, and performance of natural
estimates, we obtain
\[
\tilde F_1(\xi,\zeta,y,t)=
E_0(\xi ,\zeta ,y)\sum_{n=0}^{N-1}\sum_{j=0}^{N-n-1}\zeta^n \xi^j
{\psi_{nj}(y)\over t^{(n+j+3)/2}} (1+\delta_n(t))+\Delta_N(\xi,\zeta,y,t),
\]
where
\begin{equation}
\psi_{nj} (y)=
{ g(y)j_0(y)a^{n+j+1}(y) \over 2n!j!}\Gamma \left( {n+j+1\over 2}\right)
(1+(-1)^{n+j}),
\l{3.8}
\end{equation}
and $|\delta_n(t)|\leq \displaystyle{{B_{nj}\over \sqrt t}}$,
\begin{equation}
|\Delta_N|\leq A(N)g(y)\displaystyle{\sum_{j=0}^N{\vert\zeta^j
\xi^{N-j}\vert
\over t^{(N+3)/2}}}e^{-q_0(y)(\xi +\zeta )},
\l{3.9}
\end{equation}
where $A(N)\leq\displaystyle{\left ({e^{N+3}\over
       (N+3)^{N+3}2^{N+4}}\right )^{1/2}}$.
Estimate (\r{3.9}) is valid as $\vert\xi\vert <t^{1/4 }$,
$\vert\zeta\vert<t^{1/4}$.

Let us estimate now the second kernel $F_2(\xi,\zeta,y,t)$ containing
the integration over the set $\Omega_{\varepsilon^\prime}$.
Taking into account that  $C(y)-f(p,q,y)\geq{\varepsilon^\prime
\over 2\varepsilon}
\quad ((p,q)\in\Omega_{\varepsilon^\prime},\,\, q\geq\varepsilon >0)$
and Condition C,
we can write
\begin{align}
|F_2(\xi ,\zeta ,y,t)|
&\leq
\iint_{{\displaystyle\Omega_{\varepsilon^\prime}}}
e^{-q(\xi +\zeta )-2q(C(y)-f(p,q,y))t}\ d\mu (p,q)\notag\\
&\leq
e^{-{\varepsilon\over \varepsilon^\prime}t}
\iint_{{\displaystyle\Omega_{\varepsilon^\prime}}}
e^{q(|\xi |+|\zeta |)}\ d\mu (p,q)=O\left (e^{-{\varepsilon\over
\varepsilon^\prime}t} \right).
\l{3.91}
\end{align}
Thus assuming $\vert\xi\vert <t^{1/4 }$, $\vert\zeta\vert <
t^{1/4 }$  we obtain the final asymptotic formula
\begin{align}
  F(x ,z ,y, t)
={}&E_0(\xi ,\zeta ,y)\sum_{n=0}^{N-1}
\sum_{j=0}^{N-n-1}\zeta^n \xi^j {\psi_{nj}(y)\over t^{(n+j+3)/2}}
(1+\delta_n(t)) \notag\\
&{}+O\left (\sum_{j=0}^N{\vert\zeta^j \xi^{N-j}\vert
\over t^{(N+3)/2}}e^{-q_0(y)(\xi +\zeta )}\right )+
O\left (e^{-{\varepsilon\over
\varepsilon^\prime}t} \right)\l{392}
\end{align}
\noindent with $\xi =x-C(y)t$, $\zeta =x-C(y)t$ and
$|\delta_n(t)|\leq \displaystyle{ t^{-1/2}}$. It is not difficult to
see that this asymptotic expression can be differentiated with respect
to $x$.

The following estimates hold:
\begin{equation}
\int_x^{\infty} \int_x^{\infty} |\Delta_N(s,z,y,t)|^2dsdz \leq {
A^2(N)\over  t^{1/2-\varepsilon_1}}, \qquad (0<\varepsilon_1 <1/2)
\l{3.10}
\end{equation}
in the domain
\[
\zeta >\xi>-{1\over 2q_0(y)}\ln{ t^M\over g(y)},\qquad
|\ln g(y)|<\ln t,\qquad M={2N+5\over 4},
\]
  and
\begin{equation}
\int_x^{\infty}\int_x^{\infty}|F_2(s,z,y,t)|^2dsdz =
O(e^{-\varepsilon^\prime t}),
\l{3.11}
\end{equation}
in the domain $x>C(y)t-\sqrt t$.

Taking into account (\r{392}), (\r{3.10}) and (\r{3.11}) we can
formulate the following lemma.

\begin{lemma}
Inside the domain
\begin{equation}
\zeta >\xi >-{1\over 2q_0(y)}\ln{ t^{M+1}},
\quad|\ln g(y)|<\ln t,\quad M>1,
\l{3.111}
\end{equation}
as $t\to\infty$ the kernel $\tilde F(\xi,\zeta,y,t)$ is
represented in the form
\[
\tilde F(\xi,\zeta,y,t)=F_N(\xi,\zeta,y,t)+\tilde G(\xi,\zeta,y,t),
\]
where
\[
F_N(\xi,\zeta,y,t)=e^{ip_0(\xi-\zeta)-q_0(\xi+\zeta)}
\sum_{n=0}^{N-1}\sum_{j=0}^{N-n-1}
\zeta^n\xi^j{\psi_{nj}(y)\over t^{(n+j+3)/2}},
\]
$N=[(4M-5)/2]$, and the functions $\psi_{nj}(y)$ are bounded and defined in
$(\r{3.8})$.
The function $\tilde G(s,z,y,t)$ admits the uniform estimate
with respect to $y$ in $(\r{3.111})$:
\begin{equation}
\int_\xi^{\infty}\int_\xi^{\infty}|\tilde G(s,z,y,t)|^2dsdz=
O\left({1\over t^{1/2-\varepsilon_1}}\right)\qquad(0<\varepsilon_1 <1/2).
\l{3.12}
\end{equation}
\end{lemma}

\medskip

\noindent
\underline{Second step.}
We show that after  replacing the
kernel  $F(x,z,y,t)$ by the degenerate kernel  $F_N(x,z,y,t)$ one can obtain
an asymptotic representation of the solution
$K(x,z,y,t)$ of the Marchenko integral equation (\r{1.2})  up to
$O(t^{-1/2+\varepsilon_1})$ $(0<\varepsilon_1 <1/2)$ as $t\to\infty$.
   Set up
$\xi=x-C(y)t$, $\zeta =z-C(y)t$ and consider the domain
$z>x>C(y)-\displaystyle{{1\over 2q_0(y)}}\ln t^{M+1}$,
$|\ln g(y)|<\ln t$ with an arbitrary number $M>1$. Let
us introduce  the operators
\[
(\hat F_Nf)(z)=\int_x^\infty F_N(s,z,y,t)f(s)ds,\qquad
(\hat G_Nf)(z)=\int_x^\infty G_N(s,z,y,t)f(s)ds,
\]
in $L^2[x,\infty)$.
  Here $F_N(x,z,y,t)$ is the degenerate kernel
\begin{equation}
F_N(\xi ,\zeta ,y, t)=
e^{ip_0(\xi -\zeta )-q_0(\xi +\zeta  )}\sum_{n=0}^{N-1}
\sum_{j=0}^{N-n-1}\zeta^n \xi^j {\psi_{nj}(y)\over t^{(n+j+3)/2}},\quad
N=[(4M-5)/2],
\l{3.13}
\end{equation}
with $\xi=x-C(y)t$, $\zeta=z-C(y)t$,
and $\psi_{nj}(y)$  defined  by (\r{3.8}). $G_N(x,z,y,t)$ is the
difference between $F(x,z,y,t)$ (\r{1.4}) and $F_N(x,z,y,t)$:
\[
G_N=F-F_N.
\]
Now, (\r{1.2}) acquires the  form
\begin{equation}
(I+\hat F_N)f +\hat G_Nf =h_N +g_N,
\l{3.14}
\end{equation}
where $f=K(x,z,y,t)$, $h_N=-F_N(x,z,y,t)$, $g_N=-G_N(x,z,y,t)$.
By virtue of (\r{3.12})  one can easily  obtain the following estimates for
the norms of the operator $\hat G_N$ and the vector $g_N$ in
$L^2[x,\infty)$ as
$z >x >C(y)-\displaystyle{{1\over 2q_0(y)}}\ln t^{M+1}$, $|\ln g(y)|<\ln
t$, $t\to\infty$:
\begin{equation}
||\hat G_N||\leq A_1t^{-1/2+\varepsilon_1 },\quad
||\hat g_N||\leq A_1t^{-1/2+\varepsilon_1 }\quad
(A_1=\mbox{\rm const},\,\,\, 0<\varepsilon_1 <1/2).
\l{3.15}
\end{equation}
The operator $I+\hat F_N$ is the direct sum of the two operators
$I_1+\hat F_N$
and $I_2$. The first one acts  in the subspace $H_1$ of   $L^2[x,\infty)$,
which is generated by the vectors
\[
e^{(ip_0-q_0)z},ze^{(ip_0-q_0)z},\dots,
z^{N-1}e^{(ip_0-q_0)z}.
\]
The second operator $I_2$ $(I=I_1\oplus I_2)$  acts in the
orthogonal complement $H_2=L^2[x,\infty)\ominus H_1$.
Since we have $\hat F_N=\hat F-\hat G_N$ and the operator $I+\hat F$ is
invertible, we deduce that the
operator $I+\hat F_N$ is also invertible in $L^2[x,\infty)$   and
\begin{equation}
||(I+\hat F_N)^{-1}||_{L^2[x,\infty)}\leq A_2.
\l{3.16}
\end{equation}
We shall look for a solution of (\r{3.14}) of the form $f=\phi_N +\psi_N$,
where
$\phi_N$ is the solution of the equation $(I+\hat F_N)\phi_N=h_N$.
It implies that $\psi_N$  satisfies
\[
(I+\hat F_N)\psi_N=g_N-\hat G_N\phi_N.
\]
According to (\r{3.16}), we have
\[
||\psi_N||\leq A_2(||g_N||+||\hat G_N||\,||\psi_N||).
\]
It follows from considerations presented below  that
$\phi_N$ is uniformly bounded with respect to $(x,y)\in G_M(t)$ and $ t$
in the space $L^2[x,\infty)\cap C[x,\infty)$.
This fact and (\r{3.15}) allow us to conclude
that as $t\to\infty$
\begin{equation}
f(z)=\phi_N(z)+O(t^{-1/2+\varepsilon_1 }),\qquad 0<\varepsilon_1<1/2.
\l{3.17}
\end{equation}

\medskip

\noindent
\underline{Third step.}
The replacement of the kernel $F$  by the degenerate kernel $F_N$ in  equation
(\r{1.2}) and the implementation of the substitutions
$x=C(y)t+\xi,\,\,z=C(y)t+\zeta$ allow us to obtain the following integral
equation for the function $K_N(\xi ,\zeta , y,t)=K(\xi +C(y)t,\zeta
+C(y)t,y,t)\,\,(\zeta>\xi)$:
\begin{equation}
K_N(\xi ,\zeta , y, t)+F_N(\xi ,\zeta , y, t)+\int_\xi^\infty
K_N(\xi ,s , y, t)
F_N(s,\zeta , y, t)ds =0,\quad
\l{3.18}
\end{equation}
where $F_N$ is given by (\r{3.13}).

According to (\r{1.1}) and (\r{3.17}), the following representation of the
function $v(x,y,t)$ is valid in $G_M(t)$:
\begin{equation}
v(x,y,t)=2{\partial\over \partial\xi }
K_N(\xi ,\xi ,y,t)\bigg |_{
\xi=x-C(y)t}+O(t^{-1/2+\varepsilon_1}).
\l{3.19}
\end{equation}
We shall look for a solution of  equation (\r{3.18}) of the form
\begin{equation}
K_N(\xi ,\zeta, y, t)=\sum_{n=0}^{N-1}\gamma_n(\xi
,y,t)\zeta^ne^{-(ip_0+q_0)\zeta}.
\l{3.20}
\end{equation}
After  substitution of (\r{3.20}) into (\r{3.18}) we obtain a
system of algebraic equations for the function $\gamma_n(\xi,y,t)$:
\begin{gather}
\gamma_n+\sum_{m=0}^{N-1}\gamma_m\sum_{j=0}^{N-n-1}
{\psi_{nj}(y)\over t^{(n+j+3)/2}}\int_{\xi}^{\infty}s^{j+m}e^{-2q_0s}ds
=-\sum_{j=0}^{N-n-1}{\psi_{nj}(y)\over
t^{(n+j+3)/2}}\xi^je^{(ip_0-q_0)\xi},
\l{3.21}
\\
n=0,\dots,N-1.\notag
\end{gather}
The solution of  (\r{3.21}) has the form
\[
\gamma_l={\det [I+A(\xi ,y,t)]^{(l)}\over \det [I+A(\xi ,y,t)]},
\]
where $I$ is the identity matrix, $A(\xi,y,t)$ is the matrix with entries
\begin{equation}
\begin{split}
&[A]_{n+1,m+1}=\sum_{j=0}^{N-n-1}{\psi_{nj}\over
t^{(n+j+3)/2}}I_{j+m},\\
&I_{j+m}=\int_{\xi}^{\infty}s^{j+m}e^{-2q_0s}ds,
\quad n,m=0,\dots,N-1.
\end{split}
\l{3.22}
\end{equation}
The matrix $[I+A(\xi ,y,t)]^{(l)}$ is obtained via the substitution of the
column of right-hand sides   of the system (\r{3.21})
instead of $l$-th column of the matrix $[I+A(\xi ,y,t)]$. The functions
$\psi_{nj}$
are defined by (\r{3.8}).

The substitution of $\gamma_n(\xi ,y,t)$ into (\r{3.20}) gives us
\[
K_N(\xi ,\zeta ,y,t )=\sum_{n=0}^{N-1}
{\det (I+A)^{(n)}\over \det (I+A)}\zeta^ne^{-(ip_0+q_0)\zeta} .
\]
Hence for the solution $v(x,y,t)$ the following asymptotic representation
holds true
\begin{equation}
v(x,y,t)=2{\partial^2\over \partial x^2 }\ln \det
[I+A(x-C(y)t,y,t)]+O(t^{-{1\over 2}+\varepsilon_1})
\qquad (t\rightarrow \infty).
\l{1.52}
\end{equation}

To obtain  asymptotics (\r{1.52}) we need to investigate the behaviour of
$\Delta=\Delta(\xi,y,t)=\det [I+A(\xi,y,t)]$ in the domain
(\r{3.111}) as  $t\to\infty$. An analysis of the
structure of the matrix $I+A(\xi ,y,t)$ shows that $\Delta(\xi ,y,t)$
can be represented in the form
\begin{equation}
\Delta(\xi ,y,t)=1+\sum_{n=1}^N{P_n(\xi,y,t)\over t^{n(n+2)/2}}
e^{-2nq_0(y)\xi},
\l{53}
\end{equation}
where $P_n(\xi,y,t)$ are polynomials with respect to $\xi$
of degree at most $N^n$ with coefficients bounded with respect to
$y$ and $t$. We have the following asymptotic relations as
$n\leq [(N+1)/2]$
\begin{equation}
{P_n(\xi , y, t)\over t^{n(n+2)/2}}e^{-2nq_0(y)\xi}=
\det [A^{(n)}(\xi ,y,t)](1+O(t^{-{1\over 2}+\varepsilon_1})),
\l{54}
\end{equation}
where $A^{(n)}(\xi,y,t)$ are the $n$ by $n$ matrices with entries
\[
A^{(n)}_{i+1,k+1}=\sum_{j=0}^{n-1}{\psi_{ij} I_{j+k}\over
t^{(i+j+3)/2}}.
\]
The matrix $A^{(n)}(\xi ,y,t)$ is obviously written as the product of
two matrices. Therefore we have from (\r{3.8}) and (\r{3.22}):
\begin{equation}
\det [A^{(n)}(\xi ,y,t)]={g^n(y)j_0^n(y)a^{n^2}(y)\over
2^n(q_0(y))^{n^2}\prod_{k=0}^{n-1}(k!)^2}
{\Omega^{(n)}\tilde I^{(n)}(\xi , y)\over t^{{n(n+2)\over 2}}},
\l{55}
\end{equation}
where $\Omega^{(n)}$ and $\tilde I^{(n)}(\xi, y)$ are the determinants
of the $n$ by $n$ matrices with entries
\begin{gather}
\begin{split}
&\Omega^{(n)}_{i+1,k+1}=\Gamma\left ({i+k+1\over 2}\right )
\left (1+(-1)^{i+k}\right ),\\
&\tilde I^{(n)}_{i+1,k+1}=\tilde I_{i+k}=
\int_{q_0(y)\xi}^\infty s^{i+k}e^{-2s}ds,
\end{split}
\l{56}
\\
i,k=0,1,\dots,n-1. \notag
\end{gather}
The elements of the determinant $\tilde I^{(n)}$ satisfy the relations
\[
{d\tilde I_{i+k}\over d\xi}
=-2q_0(y)\tilde I_{i+k}+(i+k)q_0(y)\tilde I_{i+k-1},
\]
whence we have the equality
\[
{d\tilde I^{(n)}\over d\xi}+2nq_0(y)\tilde I^{(n)}=q_0(y)
\sum_{k=0}^{n-1} B^{(n)}_{k+1},
\]
where $B^{(n)}_{k+1}$ is the determinant obtained from $\tilde
I^{(n)}$ by replacing the $(k+1)$-th column by the column with
entries $i\tilde I_{i+k-1}$ $(i=0,1,\dots,n-1)$. We can show that
$\sum_{k=0}^{n-1} B^{(n)}_{k+1}\equiv 0$, hence
$\tilde I^{(n)}(\xi,y)=\tilde I^{(n)}(0 , y)e^{-2nq_0(y)\xi}$.
It follows from (\r{56}) that $\tilde I^{(n)}(0,y)$ do not depend
on $y$, therefore we write finally
\[
\tilde I^{(n)}(\xi , y)=\tilde I^{(n)}_0e^{-2nq_0(y)\xi},
\]
where $\tilde I^{(n)}_0$ are positive numbers depending on $n$.

The determinant $\tilde I^{(n)}_0$ is, up to a factor $2^{-n^2}$,
the Gram determinant of the system of functions $x^ke^{-x/2}$
$(k=0,1,\dots,n-1)$ on the semi-axis $[0,\infty)$, therefore
$\tilde I^{(n)}_0\ne 0$.  The determinants $\Omega^{(n)}$ are also
the Gram determinants of the system of functions
\[
u_k=
\begin{cases}(-1)^{(k+1)/4}x^{(2k-1)/4}e^{x/2},\quad\hfill -\infty
   <x<0\\
x^{(2k-1)/4}e^{-x/2},\quad\hfill 0<x<\infty
\end{cases}
\]
on the axis $(-\infty ,\infty)$, and $\Omega^{(n)}\ne 0$. Since in the
region (\r{3.111}) the determinant $\Delta$ does not
vanish (the operator $(I+F)^{-1}$ is bounded), it is not difficult to
see, using (\r{53})-(\r{55}), that $\Omega^{(n)}\tilde I^{(n)}_0>0$.

Thus, for $n\leq [(N+1)/2]$, $P_n(\xi ,y, t)=P_n(y, t)$ does not
depend on $\xi$, where
\begin{equation}
P_n(y,t)={g^n(y)j_0^n(y)a^{n^2}(y)\Omega^{(n)}\tilde I^{(n)}_0\over
2^n(q_0(y))^{n^2}\prod_{k=0}^{n-1}(k!)^2}
(1+O(t^{-1/2}))>0.
\l{57}
\end{equation}

Let us turn back to the representation (\r{1.52}) of the JE-I
solution. We can write
\begin{equation}
v(x,y,t)\sim 2{\partial^2\over\partial \xi^2}\ln\Delta =
2{\Delta^{\prime\prime}\Delta -(\Delta^{\prime})^2\over \Delta^2}.
\label{58}
\end{equation}
Let us cover the domain $G_M(t)$ $(M>2)$ by the subdomains
\begin{gather*}
a_1=\left\{-{1\over 2q_0}\ln{t^{2+\varepsilon}\over g(y)}<\xi
<\infty\right\},
\\
a_n=\left\{-{1\over 2q_0}\ln {t ^{n+1+\varepsilon}\over
g(y)}<\xi <-{1\over 2q_0}\ln {t^{n-\varepsilon}\over g(y)}\right\},
\quad n=2,3,\dots,m-1,
\\
a_m=\left\{-{1\over 2q_0}\ln {t^M\over g(y)}<\xi <-{1\over
2q_0}\ln {t^{m-\varepsilon }\over g(y)}\right\},
\end{gather*}
where $m=\left [\displaystyle{{N+1\over 2}}\right ]=[M-1]$, $\xi = x-C(y)t$.
With $\xi$ being inside any specific $a_n$  the numerator and
denominator in (\r{58}) have the asymptotic representations
\begin{gather}
\Delta^{\prime\prime}\Delta
-(\Delta^{\prime})^2={4q_0^2(y)e^{-4(n-1)q_0(y)\xi}\over t^{(n-1)(n+1)/2}}
{g^{2n-1}\tilde P_{n-1} \tilde P_ne^{-2q_0(y)\xi}\over
t^{(2n+1)/2}}(1+O(t^{-1/2})),
\label{59}
\\
\Delta^2={g^{2(n-1)}e^{-4(n-1)q_0(y)\xi}\over t^{(n-1)(n+1)}}
\left [\tilde P_{n-1}+{g\tilde P_ne^{-2q_0(y)\xi}\over t^{(2n+1)/2}}\right
]^2(1+O(t^{-1/2}))
\label{60}
\end{gather}
with
\[
\tilde P_n(y)={j_0^n(y)a^{n^2}(y)\Omega^{(n)}\tilde I^{(n)}_0\over
2^n(q_0(y))^{n^2}\prod_{k=0}^{n-1}(k!)^2}.
\]
Using (\r{59}) and (\r{60}), we conclude that the solution
of the JE-I possesses the following asymptotic representation
uniformly with respect to $\xi$:
\begin{align}
v(x,y,t)
&
\sim{8q_0^2(y)g\tilde P_{n-1}\tilde P_n{e^{-2q_0(y)\xi}\over t^{-(2n+1)/2}
}\over
\left (\tilde P_{n-1}+g\tilde P_n{e^{-2q_0(y)\xi}\over t^{-(2n+1)/2} }\right
)^2}{\Bigg|}_{\,\xi =x-C(y)t}\notag\\
&=
{8q_0^2(y){g\tilde P_n\over \tilde P_{n-1}}{e^{-2q_0(y)\xi}\over
     t^{(2n+1)/2}}\over
\left (1+{g\tilde P_n\over \tilde P_{n-1}}{e^{-2q_0(y)\xi}\over
t^{(2n+1)/2}}\right )^2}{\Bigg|}_{\,\xi =x-C(y)t}\\
&=
{2q_0^2(y)\over \cosh^2\left [q_0(y)\left (x-C(y)t+{1\over
2q_0(y)}\left ( \ln t^{n+1/2}-\ln g(y)-\ln{\tilde P_n(y)\over \tilde
P_{n-1}(y)}\right
)\right )\right ]}.\notag
\label{61}
\end{align}
Thus, the JE-I solution splits in $G_M(t)$ into $[M-1]$ solitons of the
form
\begin{equation}
v_n(x,y,t)={2q_0^2(y)\over \cosh^2[q_0(y)\phi_n(x,y,t)]}
\l{62}
\end{equation}
with curved lines of constant phase
\[
\phi_n(x,y,t)=x-C(y)t+{1\over
2q_0(y)}\left ( \ln t^{n+1/2}-\ln g(y)-\ln\tilde\varphi_n(y)\right )
\]
where
\begin{align*}
\tilde\varphi_n(y)
&={j_0(y)a^{2n-1}(y)\over
2(q_0(y))^{2n-1}[(n-1)!]^2}{\Omega^{(n)}\tilde I^{(n)}_0\over
\Omega^{(n-1)}\tilde I^{(n-1)}_0}\\
&={(C(y)+48(C^\prime(y))^2)^{n-1}(1+24C^{\prime\prime}(y))^{n-1/2}
Q^{(n)}\Gamma^{(n)}\over
2^{(2n+5)/2}((n-1)!)^2(C(y)+12(C^\prime(y))^2)^{(10n-3)/4}Q^{(n-1)}
\Gamma^{(n-1)}},
\end{align*}
and $\Omega^{(n)}$, $\tilde I^{(n)}_0$ are the determinants of the $n$ by $n$
matrices with entries
\begin{align*}
[\Omega^{(n)}]_{i+1,k+1}&=\Gamma\left ({i+k+1\over 2}\right )
\left (1+(-1)^{i+k}\right ),\\
[\tilde I^{(n)}_0]_{i+1,k+1}&=
\Gamma (i+k+1),\quad i,k=0,1,\dots,n-1.
\end{align*}
The theorem is proved.\hfill$\qed$

\section{Examples}

\begin{example}
Assume that
\[
C(y)={y^2\over 24}+{1\over 16}.
\]
Then
\[
p_0(y)={y\over 6},\quad q_0(y)={\sqrt{y^2+1}\over 4},
\]
and the set $\Omega$ has the form
\[
\Omega=\Bigl\{ (p,q)\,\Big|\, -\infty <p<\infty,\; 0<\varepsilon \leq
q\leq {\sqrt{36p^2+1}\over 4}\Bigr\}.
\]
Determine $g(y)$ and $d\mu$ as follows
\[
g(y)=e^{-y^2},\quad d\mu =e^{-(18p^2+2q^2-1/2)}dpdq.
\]
Then Conditions A-C are fulfilled, and according to Theorem 1
there exists a JE-I solution which splits as $t\to\infty$ in the domain
($M>2$)
\[
G_M(t)=\Bigl\{(x,y)\in\D{R}^2\,\Big|\,
  |y|<\sqrt{\ln t},\;
x>{y^2t\over 24}+{t\over 16}-{2\over \sqrt{y^2+1} }\ln{ t^{M+1}}\Bigr\}
\]
into $[M-1]$ curved solitons of the form
\begin{gather*}
v_n(x,y,t)={y^2+1\over 8\cosh^2\left [{\sqrt{y^2+1}\over
       4}\psi_n(x,y,t)\right ]},\\
\psi_n(x,y,t)=
   x-{y^2t\over 24}-{t\over
   16}+{2\over \sqrt{y^2+1} }\\
\qquad\times\left ( \ln t^{n+1/2}+y^2 -\ln
   {3^{n-1/2}(6y^2+1)^{n-1}Q^{(n)}\Gamma^{(n)}\over
     2^{5n-3/2}(2y^2+1)^{(10n-3)/4}
[(n-1)!]^2Q^{(n-1)}\Gamma^{(n-1)}}\right).
\end{gather*}
\end{example}

\begin{example}
Suppose now that
\[
C(y)=b^2\quad(b=\mbox{const}>\varepsilon>0).
\]
Then
\[
p_0(y)={y\over 12},\,\, q_0(y)=b,
\]
and the set $\Omega$ has the form
\[
\Omega =\left\{(p,q)\mid -\infty <p<\infty,\quad 0<\varepsilon \leq
q\leq b\right\}.
\]
Determine $g(y)$ and $d\mu$ as follows:
\[
g(y)=e^{-y^2},\qquad d\mu =e^{-(12p)^2}dpdq.
\]
Then Conditions A-C are fulfilled, and according to Theorem 1
there exists a JE-I solution which splits as $t\to\infty$ in the domain
($M>2$)
\[
G_M(t)=\Bigl\{(x,y)\in\D{R}^2\,\Bigm|\,|y|<\sqrt{\ln t}, \;
x>b^2t-{1\over 2b^2}\ln{ t^{M+1}}\Bigr\}
\]
into $[M-1]$ curved solitons of the form
\begin{align*}
v_n(x,y,t)
&={2b^2\over \cosh^2\left [b\psi_n(x,y,t)\right ]},\\[2mm]
\psi_n(x,y,t)&= x-b^2t\\
&\quad+{1\over 2b }
\left ( \ln t^{n+1/2}+y^2 -\ln
{Q^{(n)}\Gamma^{(n)}\over b^{3(n-1/2)}[(n-1)!]^2
Q^{(n-1)}\Gamma^{(n-1)}}\right ) .
\end{align*}
\end{example}

\begin{example}
We consider the function $C(y)$ (and also the functions
$p_0(y)$, $q_0(y)$ and the set $\Omega$) introduced in the previous
example. Such structure of $\Omega$ allows us to formulate a weaker
condition on the function $g(y)$. We suppose that
$g(y)=(1+y^{2\alpha})^{-1}$ with integer $\alpha\geq 4$ and
the measure $d\mu$
instead of (\r{1.8}) satisfies the following inequality:
\[
\iint_{\Omega}{d\mu(p,q)\over 1+(12p)^{2\alpha}}
<\infty.
\]
Then the function $v(x,y,t)$ constructed by the scheme
(\r{1.1})-(\r{1.3}) satisfies the JE-I, but it is not infinitely
differentiable. We can prove that such a JE-I solution splits as
$t\to\infty$ in the domain ($M>2$)
\[
G_M(t)=\Bigl\{(x,y)\in\D{R}^2\,\Bigm|\,
  y^{2\alpha}< t, \;
x>b^2t-{1\over 2b^2}\ln{ t^{M+1}}\Bigr\}
\]
into $[M-1]$ curved solitons of the form
\begin{align*}
v_n(x,y,t)
&={2b^2\over \cosh^2\left [b\psi_n(x,y,t)\right ]},\\[2mm]
\psi_n(x,y,t)
&= x-b^2t\\
&\quad
+{1\over 2b }\left ( \ln t^{n+1/2}+\ln (1+y^{2\alpha})
-\ln
   {Q^{(n)}\Gamma^{(n)}\over b^{3(n-1/2)}[(n-1)!]^2
Q^{(n-1)}\Gamma^{(n-1)}}\right ) .
\end{align*}
Their lines of constant phase are deviated from the straight line just
on the value ${1\over 2}\ln (1+y^{2\alpha})$. Therefore we call them weakly
curved solitons.
\end{example}

The approach developed in the present paper can be applied to solve inverse
problems for other dispersion models, like the first-order Debye model, or
more generally, the multi-resonance Lorentz and $N\,$th-order Debye models.
The generalization is straightforward and is based on the construction of
proper piecewise holomorphic functions near each pole of the model.


\label{lastpage}
\end{document}